# RATE OF RELAXATION FOR A MEAN-FIELD ZERO-RANGE PROCESS[1]

By Benjamin T. Graham

*University of British Columbia*

We study the zero-range process on the complete graph. It is a Markov chain model for a microcanonical ensemble. We prove that the process converges to a fluid limit. The fluid limit rapidly relaxes to the appropriate Gibbs distribution.

**1. Introduction.** Suppose there are a number of boxes, $N$, each containing $R$ indistinguishable balls. At rate $N$, do the following. Pick a *source* box and a *sink* box—do this uniformly at random over the $N^2$ ordered pairs of boxes. If the source box is not empty, take a ball from there and place it in the sink box. This is a Markov chain on the set

$$\mathcal{B}_N = \{\mathbf{b} \in \mathbb{N}^N : b_1 + \cdots + b_N = NR\}.$$

The number of balls in the $i$th box is $b_i$. We will write $B_i = B_i(t)$ for the corresponding random variable.

We will call this the mean-field zero-range process (MFZRP). The zero-range process is normally defined on a directed graph, with balls jumping along the edges. In contrast, our process is implicitly defined on a complete graph—balls can move from any box to any other box. The transition rate between neighboring elements of $\mathcal{B}_N$ is $N^{-1}$. The Markov chain is reversible with respect to the uniform distribution on $\mathcal{B}_N$. In the language of statistical physics, the equilibrium process can be said to have Bose–Einstein statistics [6, 9].

The number of ways of putting $b$ indistinguishable balls into $N$ distinguishable boxes is $\binom{b+N-1}{N-1}$; the number of configurations $|\mathcal{B}_N| = \binom{NR+N-1}{N-1}$.

Received December 2007; revised May 2008.
[1]Supported by an EPSRC Doctoral Training Award to the University of Cambridge.

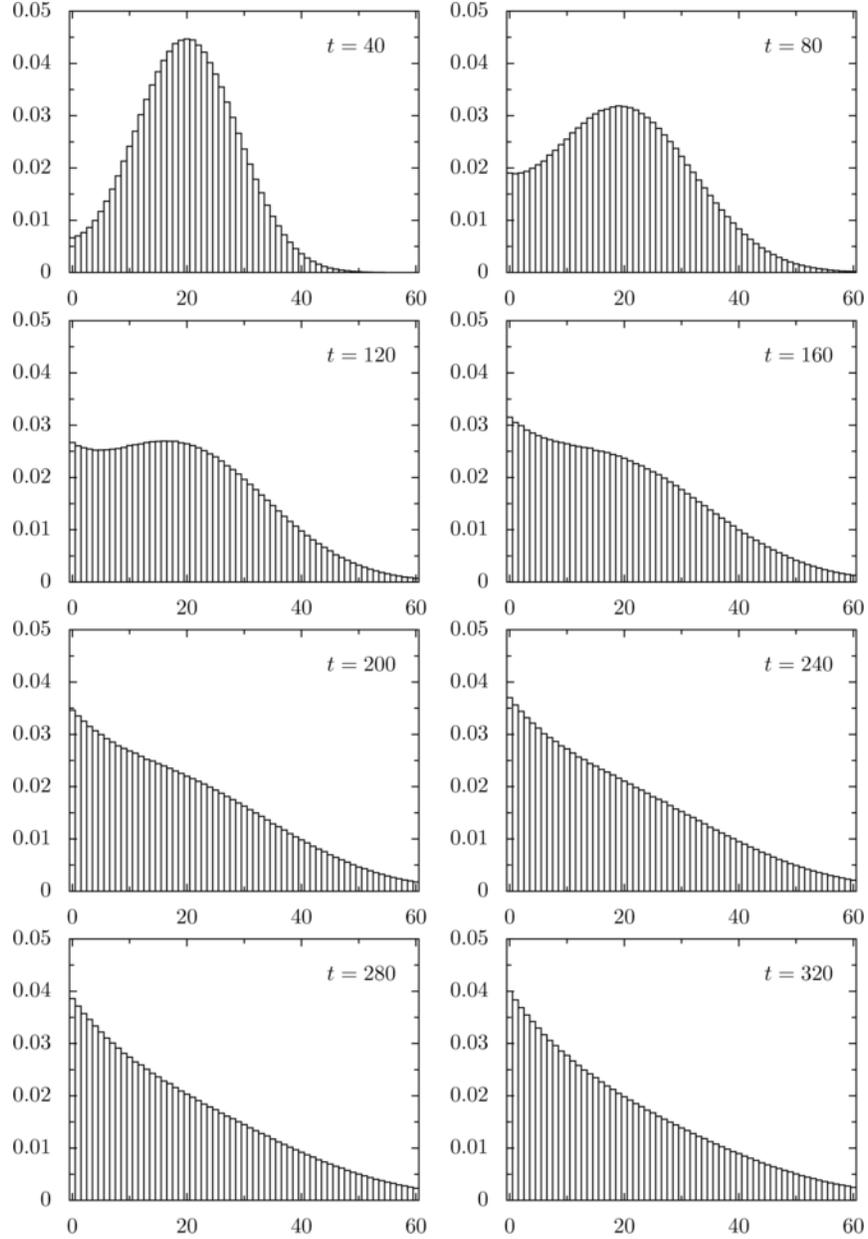

FIG. 1. *A sample path: the empirical distribution with $R=20$, $N=10^7$.*

Under the equilibrium measure $\pi$,

$$\pi(B_1 \geq k) = \binom{NR - k + N - 1}{N - 1} \bigg/ \binom{NR + N - 1}{N - 1}$$

$$= \left(\frac{R}{R+1}\right)^k (1 + \mathrm{O}(N^{-1})) \qquad \text{as } N \to \infty.$$

Let $N_k$ be the number of boxes containing $k$ balls, and let $X_k = N_k/N$. We will call $\mathbf{X} = (X_0, X_1, \ldots)$ the empirical distribution of the process. In Figure 1, we show a series of snapshots of the empirical distribution as it evolves with time. The process converges to a fluid limit as $N \to \infty$. The fluid limit methodology has been applied to a large number of models [2, 7, 15]. For an introduction to fluid limits, and a general fluid limit theorem in finite dimensions, see [2]. We will use an extension of this theorem to infinite dimensions [7], Chapter 2.

THEOREM 1.1. *Up to any finite time $T$, the empirical distribution converges exponentially, in sup norm $\|\cdot\|$, to a fluid limit $\mathbf{x}$. For $\delta > 0$,*

$$\limsup_{N \to \infty} N^{-1} \log \mathbb{P}\left(\sup_{0 \le t \le T} \|\mathbf{X}(t) - \mathbf{x}(t)\| > \delta\right) < 0.$$

We will study the MFZRP by looking at the fluid limit. We saw above that the equilibrium distribution of $B_1$ converges to a geometric distribution as $N \to \infty$. We therefore expect the fluid limit to do the same as $t \to \infty$.

THEOREM 1.2. *Let $G^R$ be the geometric distribution with mean $R$. The fluid limit converges to $G^R$ exponentially in Kullback–Leibler divergence, and*

$$\limsup_{t \to \infty} t^{-1} \log D_{\mathrm{KL}}(\mathbf{x}(t) \| G^R) = -\Omega(R^{-2}).$$

The motivation for the MFZRP comes from statistical physics. In Section 3, we introduce the microcanonical ensemble formalism [5, 12]. The fundamental assumption of statistical physics is that equivalent states have equal probability. Ensembles were introduced to derive the Boltzmann and Gibbs distributions from the assumption of equiprobability.

The MFZRP has also been studied in connection with the ZRP on the discrete torus, $\mathbb{Z}^d / L\mathbb{Z}^d$ [16]. By constructing "flows" between different graphs, the rates of relaxation of the corresponding ZRP can be compared. The spectral gap of the MFZRP is $\Omega(R^{-2})$ uniformly in $N$.

We can look at the boxes as distinguishable particles, and the balls as indistinguishable quanta of energy. At equilibrium, the probability that $(N_0, N_1, \ldots) = (n_0, n_1, \ldots)$ is proportional to the number of configurations $\mathbf{b} \in \mathcal{B}_N$ compatible with $(n_0, n_1, \ldots)$. From this point of view, the MFZRP is a microcanonical ensemble. The corresponding Gibbs distribution is the geometric distribution $G^R$.

In Section 4, we describe a simple Markov chain microcanonical ensemble, the Ehrenfest model [4, 5]. It was proposed as a probabilistic model for entropy. The Ehrenfests wanted to reconcile the deterministic nature of Boltzmann's $H$-theory with the apparently random microscopic disorder of matter.

In Section 5, we derive the MFZRP fluid limit. The differential equation that defines the fluid limit has a unique stationary point, $G^R$. Imagine that $X_0$, the fraction of empty boxes in the finite MFZRP, is "fixed." Then any given box behaves like a biased random walk on $\mathbb{N}$. It loses balls at rate 1, and gains balls at rate $1 - X_0$. This type of biased random walk is central to our analysis. In Section 8, we use properties of the biased random walk to prove bounds on the fluid limit.

An important feature of the biased random walk is that it can be truncated to a finite Markov chain without changing the "typical" behavior. In Section 6, we discuss some general convergence techniques for reversible Markov chains [3, 17]. Coercive inequalities force finite Markov chains to converge to equilibrium exponentially fast in both $\chi^2$ distance and Kullback–Leibler divergence. We consider variants of these inequalities in Section 7. These will be used in Section 9 to show that the fluid limit converges exponentially.

In Section 10, we prove the convergence of the empirical distribution to the fluid limit. The proof provides a rigorous justification for the fluid limit differential equation derived in Section 5.

Our results are stated for the MFZRP started with $R$ balls in each box, and the corresponding fluid limit as $N \to \infty$. We can extend the results to allow different initial distributions of the $NR$ balls. Suppose in the limit as $N \to \infty$, $\|\mathbf{X}(0) - \mathbf{x}(0)\| \to 0$, where $\mathbf{x}(0)$ is a probability distribution with mean $R$. The proof of Theorem 1.1 is easily adapted. However, the rate of convergence of the fluid limit to $G^R$ will depend on the tail of $\mathbf{x}(0)$; $D_{\mathrm{KL}}(\mathbf{x}\|G^R)$ will decay exponentially only if there is a geometric bound on the tail of $\mathbf{x}(0)$.

A similar dependence on the initial configuration seems to arise when considering the total variation mixing time for the finite MFZRP. The total variation mixing time $\tau_1(1/4) = \mathrm{O}(NR^2 \log R)$ [16]. However, we may expect that the process mixes much more quickly when started with $\mathbf{X}(0) = \delta_R$. We discuss this further in Section 11.

**2. Notation.** Let $\Omega$ be a countable set. For convenience, we will use the terms "measure" and "mass function" interchangeably. Let $\delta_\ell$ be the Dirac measure with its atom at $\ell$:

$$\delta_\ell(k) = \begin{cases} 1, & k = \ell, \\ 0, & \text{otherwise.} \end{cases}$$

We will say that $Q: \Omega \times \Omega \to \mathbb{R}$ is a transition rate matrix ($Q$-matrix) if for all $j \neq k$, $Q(j,k) \geq 0$ and for all $k$, $Q(k,k) = -\sum_{j \neq k} Q(k,j)$. Assume now that $\sup_k |Q(k,k)|$ is finite. The distribution at time $t$ of the continuous-time Markov chain generated by $Q$ with initial distribution $\mu$ is

$$\mu_t = \mu e^{tQ}, \qquad \frac{d}{dt}\mu_t = \mu_t Q.$$

Define $\|\mu - \pi\|_{\mathrm{TV}}$ to be the total variation distance between measures $\mu$ and $\pi$ on $\Omega$. We will also need a stronger norm to measure the distance between distributions on the nonnegative integers with finite means. We will write $\|\mu - \pi\|_{(1)}$ for the first moment of $|\mu - \pi|$,

$$\|\mu - \pi\|_{(1)} = \sum_k k|\mu(k) - \pi(k)| \geq \frac{1}{2}\sum_k |\mu(k) - \pi(k)| = \|\mu - \pi\|_{\mathrm{TV}}.$$

The Kullback–Leibler divergence from $\mu$ to $\pi$, also called the relative entropy of $\mu$ with respect to $\pi$, is defined up to a multiplicative factor by

$$D_{\mathrm{KL}}(\mu\|\pi) = \sum_k \mu(k)\log\frac{\mu(k)}{\pi(k)}.$$

We will take the function log to be the natural logarithm. Define a function $\phi\colon (0,\infty)\times[0,\infty)\to\mathbb{R}$,

$$\phi(x,y) = y\log y - x\log x - (y-x)\left[\frac{d}{dt}t\log t\lceil_{t=x}\right] = y\log\frac{y}{x} - (y-x).$$

If $x = y$, $\phi(x,y) = 0$. By the convexity of $t\mapsto t\log t$, $\phi(x,y) \geq 0$ for all $x$, $y$. The function is bounded above: $\phi(x,y) \leq (x-y)^2/x$. We can write

$$D_{\mathrm{KL}}(\mu\|\pi) = \sum_k \phi(\pi(k), \mu(k)) \geq 0.$$

Let $\mu$ and $\nu$ be probability distributions on the extended real number line $\mathbb{R}\cup\{-\infty, +\infty\}$. We will say that $\mu$ is *stochastically smaller* than $\nu$, written $\mu \leq_{\mathrm{st}} \nu$, if the following equivalent [11] conditions are met:

(i) For all $c \in \mathbb{R}$, $\mu(\cdot \geq c) \leq \nu(\cdot \geq c)$.
(ii) There is a coupling $(X,Y)$ such that $X$ has distribution $\mu$, $Y$ has distribution $\nu$, and $X \leq Y$ almost surely.

We will make use of asymptotic notation: "order less than," O; "order equal to," $\Theta$; and "order greater than," $\Omega$, to describe the asymptotic behavior of functions. With $c_1$, $c_2$ and $n_0$ positive constants, if $\forall n \geq n_0$, $|f(n)| \leq c_2|g(n)|$      write $f(n) = \mathrm{O}(g(n))$ as $n\to\infty$, if $\forall n \geq n_0$, $c_1|g(n)| \leq |f(n)| \leq c_2|g(n)|$      write $f(n) = \Theta(g(n))$ as $n\to\infty$, if $\forall n \geq n_0$, $c_1|g(n)| \leq |f(n)|$      write $f(n) = \Omega(g(n))$ as $n\to\infty$.

Likewise, $f(x) = \mathrm{O}(g(x))$ as $x\to 0$ if $f(x) = \mathrm{O}(g(x))$ as $1/x \to \infty$. We will write $f(m,n) = \Theta(g(m)h(n))$ as $n\to\infty$ if:

(i) for all $m$, $f(m,n) = \Theta(h(n))$ as $n \to \infty$,
(ii) we can take the implicit constants $c_1 = \Omega(g(m))$ and $c_2 = \mathrm{O}(g(m))$.

For example, the number of balls $NR = \Theta(N)$ as $N \to \infty$, but we can also write $NR = \Theta(NR)$ as $N \to \infty$ to indicate the dependence on $R$.

**3. Microcanonical ensembles.** The Gibbs distribution is central to the study of statistical physics in the discrete setting [5, 12]. Consider a system whose state space $\Omega$ is a countable set. Suppose that the system has an energy function, or Hamiltonian, $\mathcal{H} : \Omega \to \mathbb{R}$. The Gibbs distribution at inverse-temperature $\beta$ is defined by

$$G_\beta(k) = \frac{e^{-\beta \mathcal{H}(k)}}{Z(\beta)}, \qquad Z(\beta) = \sum_{k \in \Omega} e^{-\beta \mathcal{H}(k)}.$$

A microcanonical ensemble is a collection of $N$ copies of the system, with states, say, $b_1, \ldots, b_N \in \Omega$. Canonical ensembles were used by Maxwell, Boltzmann and Gibbs to develop the theory of thermodynamics. This has aroused some interest in Markov chain models for microcanonical ensembles, such as the Ehrenfest model [4]. In Boltzmann's $H$-theory, entropy always increases. However, the universe is apparently a time-reversible system. Reversible Markov chain models have a seemingly paradoxical property. By ergodicity, any function of the system, such as entropy, that goes up must also come down.

A microcanonical ensemble evolves with time in such a way that the total energy,

$$\mathcal{E}_{\mathrm{tot}} = \sum_{i=1}^{N} \mathcal{H}(b_i),$$

is conserved. The basic assumption of statistical mechanics is that at equilibrium, all compatible configurations are equally likely. Define $C(E_{\mathrm{tot}})$ to be the set of configurations on the ensemble compatible with total energy $E_{\mathrm{tot}}$,

$$C(E_{\mathrm{tot}}) = \{(b_1, \ldots, b_N) : \mathcal{E}_{\mathrm{tot}} = E_{\mathrm{tot}}\}.$$

In the context of probability theory, the assumption is that the ensemble is an irreducible Markov chain on $C(E_{\mathrm{tot}})$, and that the equilibrium distribution is uniform.

Let $N_k = |\{i : b_i = k\}|$ be the multiplicity in the ensemble of state $k$. With $X_k = N_k/N$, the empirical distribution is $\mathbf{X} = (X_k : k \in \Omega)$. The multinomial coefficient $N!/\prod_k N_k!$ counts the permutations of the ensemble. Under the equilibrium measure $\pi(\cdot \mid \mathcal{E}_{\mathrm{tot}} = E_{\mathrm{tot}})$,

(3.1) $\quad \pi(\mathbf{X} = \mathbf{x} \mid \mathcal{E}_{\mathrm{tot}} = E_{\mathrm{tot}}) = |C(E_{\mathrm{tot}})|^{-1} \dfrac{N!}{\prod_{k \in \Omega} n_k!}, \qquad x_k = n_k/N.$

Define $S$ to be the information-theoretic entropy of a distribution $p$ on $\Omega$,

$$S(p) = -\sum_{k \in \Omega} p(k) \log p(k).$$

Assume for now that $\Omega$ is finite. The entropy of the empirical distribution is related to the thermodynamic entropy $\log[N!/\prod_k N_k!]$. The difference between the entropy of the empirical distribution, and thermodynamic entropy divided by $N$, vanishes:

$$\left| S(\mathbf{X}) - \frac{1}{N} \log \frac{N!}{\prod_{k \in \Omega} N_k!} \right| \to 0 \qquad \text{as } N \to \infty. \tag{3.2}$$

This is simply by Stirling's approximation,

$$\log n! = n \log(n/e) + \mathrm{O}(\log n).$$

Now take inverse-temperature $\beta$ such that $G_\beta$ has energy $\overline{E}$,

$$\overline{E} = G_\beta(\mathcal{H}) = \sum_k G_\beta(k) \mathcal{H}(k).$$

For all distributions $p$ that also have energy $p(\mathcal{H}) = \overline{E}$,

$$S(G_\beta) - S(p) = D_{\mathrm{KL}}(p \| G_\beta) = \sum_k \phi(G_\beta(k), p(k)) \geq 0. \tag{3.3}$$

Therefore $G_\beta$ is the maximum entropy distribution. By (3.1) and (3.2), in the limit as $N \to \infty$, almost all of the equiprobable ensemble configurations correspond to values of the empirical distributions close to $G_\beta$. Therefore, by symmetry,

$$\lim_{N \to \infty} \pi(\{b_1 = k\} \mid \mathcal{E}_{\mathrm{tot}} = N\overline{E}) = G_\beta(k).$$

Let us return now to the MFZRP specifically. It is supported on configurations with $NR$ balls in total,

$$\mathcal{B}_N = \Big\{ \mathbf{b} \in \mathbb{N}^N : \sum b_i = NR \Big\}.$$

The MFZRP therefore *is* a microcanonical ensemble. The geometric distribution $G^R$ from Section 1 can be thought of as a Gibbs distribution on state space $\Omega = \mathbb{N}$ with respect to a linear energy function, say

$$\mathcal{H} : \mathbb{N} \to \mathbb{R}, \qquad \mathcal{H}(k) = k.$$

The "support" of the empirical distribution, the set $\{k \in \mathbb{N} : N_k > 0\}$, has size $\mathrm{O}(\sqrt{NR})$. Therefore, the limit (3.2) still holds. Increasing the temperature (decreasing the inverse-temperature $\beta$) corresponds to increasing the average number of balls per box.

**4. The Ehrenfest model.** The Ehrenfest model was introduced in [4] to demonstrate that a system could be ergodic over one time scale, yet still appear irreversible over shorter time periods. The original model has inspired a number of trivial and nontrivial variants. For convenience, we will consider a continuous-time version. Start with $N$ fair coins in a row. At rate $N$, pick a coin uniformly at random and toss it. This is a Markov chain on the set $\{H = \text{heads}, \ T = \text{tails}\}^N$. The Markov chain is time reversible with respect to the uniform distribution.

The Ehrenfest model is a microcanonical ensemble with respect to the set $\Omega = \{H, T\}$ with Hamiltonian $\mathcal{H}(H) = \mathcal{H}(T) = 0$. Let $X_H$, $X_T$ be the fractions of coins with heads, tails side up, respectively. The entropy of the system is

$$S(X_H, X_T) = -X_H \log X_H - X_T \log X_T.$$

Suppose we start with the coins all tails side up. This is a highly ordered state, $S = 0$. At equilibrium, the expected value of $X_H$ is $1/2$. The expected value of the entropy is close to the maximum entropy $S_{\max} = \log 2$.

A natural way to study the evolution of the system is to take the fluid limit. Condition on $X_H(t) = x_H$ and $X_T(t) = x_T$. $X_H$ increases by $N^{-1}$ at rate $Nx_T/2$; it decreases by $N^{-1}$ at rate $Nx_H/2$. The expected change $dX_H$ over a time period $dt$ is $((x_T - x_H)/2 + \mathrm{O}(dt))\,dt$. The fluid limit is the solution $(x_H, x_T)$ to the differential equations

$$\frac{dx_H}{dt} = \frac{x_T - x_H}{2}, \qquad \frac{dx_T}{dt} = \frac{x_H - x_T}{2},$$

with boundary conditions $x_H(0) = 0$, $x_T(0) = 1$. Up to any fixed time, the empirical distribution $(X_H, X_T)$ converges in probability to $(x_H, x_T)$,

$$x_H(t) = (1 - e^{-t})/2, \qquad x_T(t) = (1 + e^{-t})/2.$$

The entropy of the fluid limit is initially 0. Maximum entropy $S_{\max} = \log 2$ is obtained in the limit as $t \to \infty$ with $x_H = x_T = 1/2$.

The fluid limit suggests $X_H = 1/2 + \mathrm{O}(1/\sqrt{N})$ after $\mathrm{O}(\log N)$. The Markov chain is a symmetric random walk on the Hamming cube. The mixing properties have been studied in great detail [3, 10, 17]. The system does mix in time $\mathrm{O}(\log N)$.

In the case of the Ehrenfest model, the components of the ensemble are evolving independently. The fluid limit is just the probability distribution of a single coin. Each coin is a finite (time-homogeneous) Markov chain, so it is immediate that the fluid limit converges exponentially in Kullback–Leibler divergence [3]. We will see in the next section that the fluid limit of the MFZRP is also the distribution of a Markov chain, but a time-inhomogeneous one.

**5. The MFZRP fluid limit.** The empirical distribution is the Markov chain $\mathbf{X}(t) = (X_0, X_1, \ldots)$. If a box has $k$ balls at time $t$, it contributes mass $N^{-1}$ to $X_k(t)$. Unless $k = 0$, the box loses balls at rate 1. The box receives balls from each of the nonempty boxes at rate $N^{-1}$; the combined rate is $1 - X_0(t)$. Initially all boxes have $R$ balls in: $\mathbf{X}(0) = \delta_R$. The fluid limit $\mathbf{x}(t)$ is the limit in probability of the empirical distribution $\mathbf{X}(t)$ as $N \to \infty$. In this section, we describe a time-inhomogeneous Markov chain that provides an intuitive interpretation of the fluid limit.

In the limit, boxes start with $R$ balls, lose balls at rate 1 and gain them at rate $1 - x_0(t)$. The rate at which any pair of boxes interact goes to zero; it is as if each box is interacting with a "mean field" that depends on time, through $x_0(t)$, but not on the box's occupancy. Define a time-inhomogeneous Markov chain $C$ on $\mathbb{N}$. Start at $R$ at time 0: let $C(0) = R$. Stepping left, decreasing $C$ by 1, corresponds to a box losing a ball. Stepping right, increasing $C$ by 1, corresponds to a box receiving a ball. If at time $t$, $C(t) > 0$, step left at rate 1. Step right at rate $1 - \mathbb{P}(C(t) = 0)$. The distribution of $C(t)$ is exactly the fluid limit $\mathbf{x}(t)$; the rigorous justification for this comes in Section 10.

The Markov chain $C$ is controlled by $\mathbb{P}(C(t) = 0) = x_0(t)$. Key to showing that the fluid limit rapidly approaches $G^R$ is a related class of time-homogeneous Markov chains. Let $a \in [0, 1]$. Call the following walk $\mathrm{BRW}[\mathbb{N}, a]$; it is a random walk on $\mathbb{N}$ with bias $a$. Step left, unless at 0, at rate 1. Step right at rate $1 - a$. The transition rate matrix $Q = Q[\mathbb{N}, a]$ is specified by

$$(5.1) \qquad Q(j, k) = \begin{cases} 1, & k = j - 1, \\ 1 - a, & k = j + 1, \\ 0, & k \notin \{j - 1, j, j + 1\}. \end{cases}$$

Let $a \in (0, 1)$. The Markov chain $\mathrm{BRW}[\mathbb{N}, a]$ is irreducible. The stationary distribution is $\pi = \pi[\mathbb{N}, a]$, $\pi(k) = a(1 - a)^k$. Of course, $\mathrm{BRW}[\mathbb{N}, a]$ is reversible with respect to $\pi$; for all $j, k$, $\pi(j) Q(j, k) = \pi(k) Q(k, j)$.

Extending our notation, one can describe the inhomogeneous Markov chain $C$ as $\mathrm{BRW}[\mathbb{N}, x_0(t)]$,

$$(5.2) \qquad \frac{d}{dt}\mathbf{x}(t) = \mathbf{x}(t) Q[\mathbb{N}, x_0(t)], \qquad \mathbf{x}(0) = \delta_R.$$

The entropy of the fluid limit $\mathbf{x} = \mathbf{x}(t)$ is

$$S(\mathbf{x}) = -\sum_{k \geq 0} x_k \log x_k, \qquad 0 \log 0 = 0.$$

The distribution associated with energy $R$ with maximum entropy, the Gibbs distribution, is the geometric distribution $G^R$. By (3.3), the Kullback–Leibler divergence from $\mathbf{x}$ to $G^R$ is the difference between the entropy and the maximum entropy,

$$D_{\mathrm{KL}}(\mathbf{x} \| G^R) = S(G^R) - S(\mathbf{x}).$$

LEMMA 5.1. *The unique stable point with mean $R$ of the differential equation (5.2) is $\mathbf{x} = G^R$. For $t > 0$, the rate of increase of the entropy is*

$$\frac{d}{dt}S(\mathbf{x}) = \sum_{k \geq 0}(x_k(1 - x_0) - x_{k+1})\log \frac{x_k(1 - x_0)}{x_{k+1}}$$

(5.3)

$$\geq \sum_{k \geq 0} \frac{(x_k(1 - x_0) - x_{k+1})^2}{\max\{x_k(1 - x_0), x_{k+1}\}}.$$

As physical intuition demands, the entropy is increasing. Further, $S(\mathbf{x})$ is strictly increasing unless $\mathbf{x} = G^R$. Showing that the fluid limit converges to $G^R$ in Kullback–Leibler divergence is equivalent to showing that $S(\mathbf{x})$ increases to $S(G^R)$.

To prove that $D_{\mathrm{KL}} = D_{\mathrm{KL}}(\mathbf{x} \| G^R)$ decays exponentially, we must show that $D_{\mathrm{KL}}$ and $dS/dt = -d(D_{\mathrm{KL}})/dt$ have more or less the same order. We can measure the Kullback–Leibler divergence from $\mathbf{x}$ to geometric distributions other than $G^R$.

LEMMA 5.2. *Let $a \in (0, 1)$. The Kullback–Leibler divergence from $\mathbf{x}$ to $\pi[\mathbb{N}, a]$ is minimized by $\pi[\mathbb{N}, a] = G^R$, that is, when $a = 1/(R + 1)$.*

Taking $a = x_0$, $D_{\mathrm{KL}}(\mathbf{x} \| \pi[\mathbb{N}, x_0]) \geq D_{\mathrm{KL}}(\mathbf{x} \| G^R)$:

(5.4) $$\sum_{k \geq 1} x_k \log \frac{x_k}{x_0(1 - x_0)^k} = \sum_{k \geq 1} \phi(x_0(1 - x_0)^k, x_k) \geq D_{\mathrm{KL}}(\mathbf{x} \| G^R).$$

In order to compare the summations in (5.3) with the summations in (5.4), we need to consider a truncated, finite version of the biased random walk. We also need to obtain certain bounds on the fluid limit.

PROOF OF LEMMA 5.1. Define $m_k = x_k(1 - x_0) - x_{k+1}$ to be the *flow* from $x_k$ to $x_{k+1}$. For convenience let $m_{-1} = 0$; no boxes have $-1$ balls. The $k$th coordinate of the fluid limit increases as mass flows from $x_{k-1}$ to $x_k$, and decreases as mass flows from $x_k$ to $x_{k+1}$,

$$\frac{d}{dt}x_k = m_{k-1} - m_k, \qquad k \in \mathbb{N}.$$

Note that the mass, 1, and mean, $R$, are conserved,

$$\frac{d}{dt}\sum x_k = 0, \qquad \frac{d}{dt}\sum kx_k = 0.$$

If $d\mathbf{x}/dt = 0$, then $m_k = 0$ for all $k$; $\mathbf{x}$ is a geometric distribution.

Differentiating the entropy,

$$\frac{d}{dt}S = -\sum_{k\geq 0}(1+\log x_k)\frac{d}{dt}x_k = \sum_{k\geq 0}m_k \log \frac{x_k}{x_{k+1}}.$$

The first line in (5.3) follows as $\sum_{k\geq 0} m_k = 0$. The lower bound on $dS/dt$ follows using

$$\frac{\log a - \log b}{a-b} \geq \frac{1}{\max\{a,b\}}, \qquad a,b > 0.$$

Apply the mean value theorem to $x \mapsto \log(x)$ on the interval from $a$ to $b$. □

PROOF OF LEMMA 5.2. If $x_k = a(1-a)^k$ for all $k$, then $\mathbf{x}$ has mean $(1-a)a^{-1}$ and entropy $S(\mathbf{x}) = -\log a - (1-a)a^{-1}\log(1-a)$. By definition,

$$D_{\mathrm{KL}}(\mathbf{x}\|G^R) = \sum x_k \log x_k/G_k^R, \qquad G_k^R = R^k/(R+1)^{k+1}.$$

Using $\sum x_k = 1$ and $\sum k x_k = R$,

$$D_{\mathrm{KL}}(\mathbf{x}\|G^R) = \log\left[\frac{a(1-a)^R}{1/(R+1)(R/(R+1))^R}\right] + D_{\mathrm{KL}}(\mathbf{x}\|\pi[\mathbb{N},a]).$$

For $a \in [0,1]$, $a(1-a)^R$ is maximized by $a = 1/(R+1)$. □

**6. Markov chain convergence.** Let $\Omega$ be a finite set. Let $Q$ be the transition rate matrix for an irreducible Markov chain, and assume that the Markov chain is time reversible with respect to equilibrium probability measure $\pi$. It is standard [17] that the eigenvectors of $Q$, suitably normalized, can be turned into an orthonormal basis for $\mathbb{R}^\Omega$. We will work with the left eigenvectors; they are orthogonal with respect to inner product

$$\langle f,g\rangle_\pi = \sum_{k\in\Omega} f(k)g(k)/\pi(k).$$

Label the eigenvalues $0 = \lambda_1 > \lambda_2 \geq \lambda_3 \geq \cdots \geq \lambda_n$ with corresponding normalized eigenvectors $F_1, \ldots, F_n$. The first eigenvector $F_1 = \pi$, the stationary distribution. For a probability measure $\mu$ on $\Omega$,

$$\mu(k) = \pi(k) + \sum_{i\geq 2} \alpha_i F_i(k), \qquad \alpha_i = \langle \mu, F_i\rangle_\pi.$$

The $\chi^2$ distance from $\mu$ to $\pi$ can be written in terms of these coefficients,

$$\sum_{k\in\Omega} \frac{[\mu(k)-\pi(k)]^2}{\pi(k)} = \langle \mu-\pi, \mu-\pi\rangle_\pi = \sum_{i\geq 2} \alpha_i^2.$$

The quantity $|\lambda_2|$ is called the spectral gap. Let $\mu_t = \mu e^{tQ}$. The eigenvalue of $F_i$ with respect to $e^{tQ}$ is $e^{t\lambda_i}$, so

$$\frac{d}{dt} \sum_k \frac{[\mu_t(k) - \pi(k)]^2}{\pi(k)} \leq 2\lambda_2 \sum_k \frac{[\mu_t(k) - \pi(k)]^2}{\pi(k)}. \tag{6.1}$$

The uniform time to convergence within $\varepsilon$ in total variation is

$$\tau_1(\varepsilon) = \inf\{t : \forall \mu, \|\mu_t - \pi\|_{\mathrm{TV}} \leq \varepsilon\}.$$

Inequality (6.1) can be used to bound the total variation convergence time,

$$|\delta_j e^{tQ}(k) - \pi(k)| \leq e^{t\lambda_2} \sqrt{\pi(k)/\pi(j)}, \tag{6.2}$$

$$\|\mu_t - \pi\|_{\mathrm{TV}} \leq e^{t\lambda_2}/\pi_{\min}, \qquad \pi_{\min} := \min_x \pi(x). \tag{6.3}$$

Log Sobolev inequalities for finite Markov chains are described in [3]. Define the Dirichlet form, and the Laplacian for $f, g : \Omega \to \mathbb{R}$,

$$\mathcal{E}(f, g) = \frac{1}{2} \sum_{j,k} \pi(j) Q(j,k) [f(j) - f(k)][g(j) - g(k)],$$

$$\mathcal{L}(f) = \sum_k \pi(k) f(k)^2 \log\left(\frac{f(k)^2}{\|f\|_{2,\pi}^2}\right).$$

The log Sobolev constant is defined

$$\alpha = \min_f \left\{ \frac{\mathcal{E}(f, f)}{\mathcal{L}(f)} : \mathcal{L}(f) \neq 0 \right\}.$$

The log Sobolev constant can be used to show convergence in Kullback–Leibler divergence,

$$\frac{d}{dt} D_{\mathrm{KL}}(\mu_t \| \pi) = -\mathcal{E}\left(\frac{\mu_t}{\pi}, \log \frac{\mu_t}{\pi}\right) \leq -4\alpha D_{\mathrm{KL}}(\mu_t \| \pi). \tag{6.4}$$

Corollary A.4 of [3] gives a bound for the log Sobolev constant in terms of the spectral gap and $\pi_{\min}$,

$$\alpha \geq \frac{(1 - 2\pi_{\min})|\lambda_2|}{\log(1/\pi_{\min} - 1)}. \tag{6.5}$$

**7. Truncating the biased random walk.** We will now consider versions of the biased random walk on sets other than $\mathbb{N}$. For $n \geq 2$, let $\mathrm{BRW}[n,a]$ be the biased random walk restricted to the set $\{0, 1, \ldots, n-1\}$. The transition rate matrix $Q[n,a]$ is still specified by (5.1). Unless at the left boundary, 0, step left at rate 1. Unless at the right boundary, $n-1$, step right at rate $1-a$. If $a \in (0,1)$, the stationary distribution $\pi = \pi[n,a]$ is given by $\pi(k) = a(1-a)^k/(1-(1-a)^n)$. The stationary distribution is uniform if $a = 0$.

Now let $a \in (0,1)$. Consider the inequalities for finite Markov chains from Section 6. Inequality (6.4), containing the log Sobolev constant, becomes

$$(7.1) \quad \sum_{k=0}^{n-2}[\mu(k)(1-a) - \mu(k+1)]\log\frac{\mu(k)(1-a)}{\mu(k+1)} \geq 4\alpha \sum_{k=0}^{n-1} \mu(k)\log\frac{\mu(k)}{\pi(k)}.$$

The spectral gap inequality (6.1) becomes

$$\sum_{k=0}^{n-2}\frac{[\mu(k)(1-a) - \mu(k+1)]^2}{\pi(k+1)} \geq |\lambda_2| \sum_{k=0}^{n-1}\frac{[\mu(k) - \pi(k)]^2}{\pi(k)}.$$

We see in the next section, Lemma 8.4, that $|\lambda_2| \geq a^2/4$. The bound (6.5) on the log Sobolev constant, Corollary A.4 of [3], gives

$$\alpha \geq \frac{(1 - 2\pi(n-1))|\lambda_2|}{\log(1/\pi(n-1) - 1)} = \Theta(a/n) \quad \text{as } n \to \infty.$$

To compare (5.3) with (5.4), we would like to set $\mu(k) = x_k$ in the above inequalities. However, $\mu$ must be a probability distribution; the restriction of $\mathbf{x}$ to the set $\{0, 1, \ldots, n-1\}$ is only a probability distribution if $n > R$ and $t = 0$. We can adapt these inequalities to cope.

LEMMA 7.1. *Let $\mu$ be a probability distribution with support $\mathbb{N}$. Let $\alpha$ be the log Sobolev constant of* BRW$[n, a]$,

$$\sum_{k=0}^{n-2}[\mu(k)(1-a) - \mu(k+1)]\log\frac{\mu(k)(1-a)}{\mu(k+1)}$$
$$\geq 4\alpha\left[\sum_{k=0}^{n-1}\mu(k)\log\frac{\mu(k)}{a(1-a)^k} + \log(1-(1-a)^n)\right].$$

LEMMA 7.2. *Let $\mu$ be a probability distribution on $\mathbb{N}$. If $\mu(0) = a$,*

$$\sum_{k=0}^{n-2}\frac{[\mu(k)(1-a) - \mu(k+1)]^2}{a(1-a)^{k+1}} \geq \frac{a^2}{4}\sum_{k=1}^{n-1}\frac{[\mu(k) - a(1-a)^k]^2}{a(1-a)^k}.$$

In Section 9, we will use Lemma 7.1 to see that $D_{\mathrm{KL}}$ does indeed become small fairly quickly. However, as $D_{\mathrm{KL}}$ decreases we must increase $n$, so our bound on $\alpha$ tends to zero. The lower bound on the log Sobolev constant is too weak to allow Lemma 7.1 to be used to prove Theorem 1.2. We complete the proof of Theorem 1.2 using Lemma 7.2. Once $D_{\mathrm{KL}}$ is small, the inequality can be used to show that $\log D_{\mathrm{KL}}$ decreases with rate order $R^{-2}$.

The lower bound on the log Sobolev constant for BRW$[n,a]$ might seem pessimistic, but it does in fact have the correct order. Let $\mu$ be the uniform distribution on $\{0, 1, \ldots, n-1\}$, and start BRW$[n,a]$ with initial distribution $\mu$. By inequality (6.4),

$$\alpha \leq -\frac{1}{4\,D_{\mathrm{KL}}(\mu_t\|\pi)}\frac{d}{dt}D_{\mathrm{KL}}(\mu_t\|\pi), \qquad \mu_t = \mu e^{tQ[n,a]},\ \pi = \pi[n,a].$$

We can apply the above inequality at $t=0$ in the limit as $n \to \infty$. Written in terms of the entropy and mean of $\mu_t$, the divergence $D_{\mathrm{KL}}(\mu_t\|\pi)$ is

$$-S(\mu_t) - \log\frac{a}{1-(1-a)^n} - \log(1-a)\sum k\mu_t(k).$$

At $t=0$, these three terms are $\Theta(\log n)$, $\Theta(\log a)$ and $\Theta(n\log(1-a))$, respectively: the divergence is $\Theta(n\log(1-a))$. Now consider the rate of change of the divergence. The entropy $S$ is maximized by the uniform distribution, so $\frac{d}{dt}S(\mu_t) = 0$ at $t=0$. The second term is constant. The semigroup $e^{tQ}$ is pulling the distribution toward zero at net rate $1-(1-a)=a$; the rate of change of the third term is $-\Theta(a\log(1-a))$. Therefore $\alpha = \Theta(a/n)$ as $n \to \infty$.

PROOF OF LEMMA 7.1.  The support of $\mu$ is $\mathbb{N}$, so $\sum_{j<n}\mu(j) > 0$. Let $\nu(k) = \mu(k)/\sum_{j<n}\mu(j)$ for $k < n$. Substitute $\nu$ into inequality (7.1).  $\square$

PROOF OF LEMMA 7.2.  We will make use of a discrete version of Hardy's inequality [14]. Let $u$ and $v$ be positive functions on $\mathbb{N}$. The inequality states that

$$\sum_{j=0}^{\infty} v(j)f(j)^2 \geq \frac{1}{4B}\sum_{j=0}^{\infty} u(j)\left(\sum_{k=0}^{j} f(k)\right)^2 \qquad \forall f \in \ell^2(v),$$

where

$$B := \sup_{k \geq 0}\left(\sum_{j=0}^{k}\frac{1}{v(j)}\right)\left(\sum_{j=k}^{\infty} u(j)\right).$$

Let $u(k) = v(k) = a(1-a)^{k+1}$ for $k \in \mathbb{N}$; this gives $B = a^{-2}$. The result follows by taking

$$f(k) = \frac{\mu(k)}{a(1-a)^k} - \frac{\mu(k+1)}{a(1-a)^{k+1}}, \qquad k = 0, 1, \ldots, n-2$$

and $f(k) = 0$ for $k \geq n-1$.  $\square$

It is also helpful to consider the effect of truncating BRW$[\mathbb{N},a]$ from the other side. Let BRW$[\mathbb{N}+k,a]$ be the walk with bias $a$ on the set $\{k, k+1, k+$ $2,\ldots\}$. The equilibrium distribution, $\pi[\mathbb{N}+k,a]$, is simply the equilibrium distribution of $\mathrm{BRW}[\mathbb{N},a]$ shifted $k$ to the right,

$$\pi[\mathbb{N}+k,a](j) = \pi[\mathbb{N},a](j-k), \qquad j = k, k+1, \ldots.$$

This provides a very simple stochastic bound on $\mathrm{BRW}[\mathbb{N},a]$.

LEMMA 7.3. $\mathrm{BRW}[\mathbb{N},a]$ *started at $k$ is stochastically smaller than the equilibrium distribution of* $\mathrm{BRW}[\mathbb{N}+k,a]$; $\delta_k e^{tQ[\mathbb{N},a]} \leq_{\mathrm{st}} \pi[\mathbb{N}+k,a]$.

PROOF. This follows by coupling. Construct walks $X \sim \mathrm{BRW}[\mathbb{N},a]$ and $Y \sim \mathrm{BRW}[\mathbb{N}+k,a]$ on the same probability space as follows. Let $X(0) = k$ and choose $Y(0)$ according to the distribution $\pi[\mathbb{N}+k,a]$.

Now introduce two Poisson processes, one with rate 1 and one with rate $1-a$, to run $X$ and $Y$ for $t > 0$. The rate 1 Poisson process corresponds to steps to the left. With each arrival of the rate 1 process:

(i) decrement $X$ by 1 unless $X = 0$,
(ii) decrement $Y$ by 1 unless $Y = k$.

If $X \leq Y$ immediately before an arrival, then $X \leq Y$ after the arrival. The rate $1-a$ Poisson process corresponds to steps to the right. With each arrival of the rate $1-a$ process, increment both $X$ and $Y$ by 1. This also preserves $X \leq Y$. At $t = 0$, $X \leq Y$ so $\mathbb{P}(\forall t, X(t) \leq Y(t)) = 1$. □

**8. Bounds on the fluid limit.** To apply Lemmas 7.1 and 7.2, we need to bound certain functions of the fluid limit. The class of Markov chains $\mathrm{BRW}[\mathbb{N},a]$ has a stochastic ordering property. Let $a, b \in [0,1]$. Define Markov chains $C_a$ and $C_b$ that evolve according to $\mathrm{BRW}[\mathbb{N},a]$ and $\mathrm{BRW}[\mathbb{N},b]$, respectively. Let $\mu_a$ be the initial distribution of $C_a$ and let $\mu_b$ be the initial distribution of $C_b$.

LEMMA 8.1. *If $a \leq b$ and $\mu_a \geq_{\mathrm{st}} \mu_b$, then there is a coupling $(C_a, C_b)$ such that $C_a \geq C_b$ with probability 1.*

Fix $s > 0$. Suppose $x_0(t) \in [a,b]$ for $t \geq s$. The distribution of the time-inhomogeneous Markov chain $\mathrm{BRW}[\mathbb{N}, x_0(t)]$, the fluid limit $\mathbf{x}$, can be compared to the distributions of the time-homogeneous Markov chains. Start $C_a$ and $C_b$ at time $s$, both with initial distribution $\mathbf{x}(s)$. Let $\mu_a$ and $\mu_b$ be the distributions of $C_a$ and $C_b$,

$$\mu_a = \mu_{a,t} = \mathbf{x}(s)e^{(t-s)Q[\mathbb{N},a]}, \qquad \mu_b = \mu_{b,t} = \mathbf{x}(s)e^{(t-s)Q[\mathbb{N},b]}.$$

LEMMA 8.2. *For $t \geq s$, $\mu_b \leq_{\mathrm{st}} \mathbf{x}(t) \leq_{\mathrm{st}} \mu_a$. If $1/(5R) \leq a \leq 1/(R+1)$,*

(8.1) $\quad -\log \|\mu_{a,t} - \pi[\mathbb{N}, a]\|_{(1)} = \Omega(R^{-2}(t-s)) \qquad$ *as $(t-s) \to \infty$.*

*Similarly for $\mu_b$ if $1/(R+1) \leq b \leq 4/5$.*

Recall that $\|\cdot\|_{\mathrm{TV}} \leq \|\cdot\|_{(1)}$. Lemma 8.2 allows us to bound

(8.2) $\begin{aligned} x_k &\leq \mu_a(C_a \geq k) - \mu_b(C_b \geq k+1) \\ &\leq (1-a)^k - (1-b)^{k+1} + \|\mu_a - \pi[\mathbb{N}, a]\|_{\mathrm{TV}} + \|\mu_b - \pi[\mathbb{N}, b]\|_{\mathrm{TV}}. \end{aligned}$

We will apply Lemma 8.2 iteratively, with bounds $a$ and $b$ improving as $s \to \infty$. As $\mathbf{x}$ tends toward $G^R$ in Kullback–Leibler divergence, $x_0$ tends toward $G^R(0)$. By Pinsker's inequality [1],

(8.3) $\qquad (x_0 - G^R(0))^2 \leq \|\mathbf{x} - G^R\|_1^2 \leq 2D_{\mathrm{KL}}(\mathbf{x}\|G^R).$

Initially, however, we will take $a$ and $b$ as follows and $s = s_0$.

LEMMA 8.3. *For $t \geq s_0 = \mathrm{O}(R^2 \log R)$:*

  (i) $x_0(t) \geq a = 1/(5R)$,
  (ii) $x_0(t) \leq b = 1 - (R+1)^{-\mathrm{O}(1)}$,
  (iii) $\mathbf{x}(t) \leq_{\mathrm{st}} \pi[\mathbb{N} + j, 1/(5R)]$ *with $j = \mathrm{O}(R \log R)$.*

Now to prove the above. The simplest way to calculate (8.1) seems to be via the truncated version of the biased random walk.

LEMMA 8.4. *For $a \in (0,1)$, $\mathrm{BRW}[n,a]$ has spectral gap $|\lambda_2| \geq a^2/4$. For all $k = 0, 1, 2, \ldots, n-1$,*

$$\|\delta_k e^{tQ[n,a]} - \pi[n,a]\|_{(1)} \leq \frac{4e^{-ta^2/4}}{a^2(1-a)^{k/2}}.$$

PROOF OF LEMMA 8.1. As $\mu_a \geq_{\mathrm{st}} \mu_b$, we can choose $C_a(0) \geq C_b(0)$. The property $C_a \geq C_b$ is preserved if we run the two Markov chains as follows. At rate 1, decrement both $C_a$ (if $C_a > 0$) and $C_b$ (if $C_b > 0$). At rate $1-b$, increment both $C_a$ and $C_b$. At rate $b-a$, increment only $C_a$. □

PROOF OF LEMMA 8.2. The stochastic bounds follow as in Lemma 8.1; the fluid limit $\mathbf{x}$ evolves according to $\mathrm{BRW}[\mathbb{N}, x_0(t)]$.

The $\|\cdot\|_{(1)}$ mixing-time for $\mathrm{BRW}[\mathbb{N}, a]$ started at $\mu$, the least $t$ such that $\|\mu e^{tQ} - \pi[\mathbb{N}, a]\|_{(1)} \leq \varepsilon$, very much depends on $\mu$. The inequality in Lemma 8.4 is uniform in $n$. We can let $n$ tend to infinity,

$$\|\delta_k e^{tQ[\mathbb{N}, a]} - \pi[\mathbb{N}, a]\|_{(1)} \leq \frac{4e^{-ta^2/4}}{a^2(1-a)^{k/2}} \to 0 \qquad \text{as } t \to \infty.$$

However, with $t$ fixed, the above bound grows exponentially in $k$. When $k$ is large, it is better to use Lemma 7.3,

$$\|\delta_k e^{tQ[\mathbb{N},a]} - \pi[\mathbb{N},a]\|_{(1)} \leq (k+a^{-1}) + a^{-1}.$$

By linearity and a triangle inequality,

$$\|\mu_a - \pi[\mathbb{N},a]\|_{(1)} \leq \sum_k x_k(s) \min\left\{k + 2a^{-1}, \frac{4e^{-(t-s)a^2/4}}{a^2(1-a)^{k/2}}\right\}.$$

The bound (8.1) now follows by applying the stochastic upper bound on $\mathbf{x}(s)$, Lemma 8.3(iii). □

PROOF OF LEMMA 8.4. Let $F_1, \ldots, F_n$ be the left eigenvalues of $Q$, normalized as in Section 6. The first eigenvector $F_1$ is the equilibrium distribution $\pi$; the corresponding eigenvalue is 0. The other $n-1$ eigenvectors of the Markov chain can be written as follows. Let $A_j = \sqrt{1-a}\exp(i\pi(j-1)/n)$ for $j = 2, \ldots, n$. Then with $c_j$ a normalizing constant, $F_j(k) = c_j \operatorname{Im}[(1-A_j)A_j^k]$. Check

$$(F_j Q)(k) = -F_j(k)|1 - A_j|^2, \qquad j = 2, \ldots, n.$$

Therefore for $j \geq 2$, the $j$th eigenvalue is

$$\lambda_j = -|1 - A_j|^2 \leq -a^2/4.$$

By (6.2),

$$|\delta_k e^{tQ}(\ell) - \pi(\ell)| \leq e^{-ta^2/4}(1-a)^{(\ell-k)/2}.$$

Multiply by $\ell$, and sum over $\ell = 0, 1, \ldots, n-1$. □

PROOF OF LEMMA 8.3(i). Let $n = 2R$. We can define an evolving probability distribution,

$$\mathbf{y} = (y_k(t)) \qquad \text{on } \{0, 1, 2, \ldots, n-1\} \cup \{+\infty\}$$

with the following properties:

(i) $\mathbf{x}(t) \leq_{\text{st}} \mathbf{y}(t)$ for all $t \geq 0$, so $x_0(t) \geq y_0(t)$,
(ii) $y_0(t) \geq 1/(5R)$ for $t \geq s_0 = \mathrm{O}(R^2 \log R)$.

Apply Markov's inequality to the fluid limit,

$$\sum_{k=0}^{\infty} k x_k = R \quad \text{so} \quad \sum_{k=0}^{n-1} x_k \geq 1/2.$$

Define probability distribution $\mathbf{y}(0)$ by $y_{n-1}(0) = 1/2$ and $y_\infty(0) = 1/2$. Let $\hat{\mathbf{y}}$ be the restriction of $\mathbf{y}$ to $\{0, 1, \ldots, n-1\}$: $\hat{\mathbf{y}} = (y_k)_{k=0}^{n-1}$. Define $\mathbf{y}(t)$ as follows:

(i) For all $t$, let $y_\infty(t) = 1/2$.
(ii) Let $Q$ be the transition rate matrix for $\mathrm{BRW}[n, 0]$; let $\hat{\mathbf{y}}(t) = \hat{\mathbf{y}}(0)e^{tQ}$.

We have $\mathbf{x}(0) \leq_{\mathrm{st}} \mathbf{y}(0)$ by the choice of $\mathbf{y}(0)$. If $\mathbf{y}$ was being acted on by $\mathrm{BRW}[\mathbb{N}, 0]$, $\mathbf{x}(t) \leq_{\mathrm{st}} \mathbf{y}(t)$ would follow by a proof similar to that of Lemma 8.1. However, $\mathbf{y}$ is being acted on by $\mathrm{BRW}[n, 0]$. We must show that as mass of $\mathbf{x}$ at $n-1$ (that is coupled to mass of $\mathbf{y}$ at $n-1$) flows to $n$, we can modify the coupling suitably. The inequality $\sum_{k<n} x_k \geq 1/2$ is preserved: for any mass of $\mathbf{x}$ (coupled to mass of $\mathbf{y}$ at $n-1$) flowing to $n$, there must be some mass of $\mathbf{x}$ below $n$ that is coupled to mass of $\mathbf{y}$ at $+\infty$. Using this slack, an exchange can be made in the coupling.

The equilibrium distribution of $\mathrm{BRW}[n, 0]$ is uniform. The mass of $\hat{\mathbf{y}}$ is $1/2$, so $y_0 \to (1/2)n^{-1} = 1/(4R)$. To see that $y_0 \geq 1/(5R)$ after time $O(R^2 \log R)$, apply (6.3). The spectral gap of $\mathrm{BRW}[n, 0]$ is given by setting $a = 0$ in the proof of Lemma 8.4,

$$|\lambda_2| = |1 - A_2|^2 \geq 4n^{-2}, \qquad A_j = \exp(i\pi(j-1)/n). \qquad \square$$

To show Lemma 8.3(iii), we must bound how far $\mathbf{x}$ shifts to the right before the lower bound $x_0 \geq 1/(5R)$ is in effect.

LEMMA 8.5.  *Let $X = Y - Z$ with $Y, Z$ independent $\mathrm{Poisson}(s_0)$ random variables. Let $\nu$ be given by $\nu(k) = \mathbb{P}(X \in \{k-R, k+R+1\})$. Then $\mathbf{x}(s_0) \leq_{\mathrm{st}} \nu$.*

We will need to apply a concentration bound to $\nu$.

THEOREM 8.6 ([13], Theorem 2.7).  *Let $X$ be the sum of $n$ independent random variables: $X = X(1) + \cdots + X(n)$. Let $\mathrm{Var}(X)$ be the variance of $X$. Suppose $X(i) - \mathbb{E}(X(i)) \leq 1$ for each $i$. For any $\lambda \geq 0$,*

$$\mathbb{P}(X \geq \mathbb{E}(X) + \lambda) \leq \exp\left(\frac{-\lambda^2}{2\mathrm{Var}(X) + 2\lambda/3}\right).$$

PROOF OF LEMMA 8.3(iii).  If $\mathbf{x}(s_0) \leq_{\mathrm{st}} \pi[\mathbb{N} + k, 1/(5R)]$, then for all $t \geq s_0$, $\mathbf{x}(t) \leq_{\mathrm{st}} \pi[\mathbb{N} + k, 1/(5R)]$. We can use the concentration bound, Theorem 8.6, to show $\nu \leq_{\mathrm{st}} \pi[\mathbb{N} + j, 1/(5R)]$ with $j = O(R \log R)$. The Poisson distributions $Y$ and $Z$ can be approximated by binomial distributions (which in turn can be written as the sums of Bernoulli random variables).  $\square$

PROOF OF LEMMA 8.5.  By Lemma 8.2 with $a = 0$, $\mathbf{x}$ is stochastically smaller than the walk $\mathrm{BRW}[\mathbb{N}, 0]$ started at $R$. We can use a reflection principle to calculate the distribution $\nu = \delta_R e^{s_0 Q[\mathbb{N}, 0]}$. Let $C$ be a random walk on $\mathbb{Z}$ that starts at $R$, then steps forward at rate 1, and backward at rate 1; in keeping with our notation, the walk is BRW$[\mathbb{Z}, 0]$. Consider the reflection of $C$ in the point $-1/2$: $-1 - C$ start at $-1 - R$ and then also evolves according to BRW$[\mathbb{Z}, 0]$. The law of $\max\{C, -1 - C\}$ is exactly the distribution of the random walk BRW$[\mathbb{N}, 0]$ started at $R$. □

PROOF OF LEMMA 8.3(ii). This follows from part (iii). Choose $n$ and $b$ such that
$$\sum_{k \geq n} k\pi[\mathbb{N}+j, a](k) \leq R, \qquad b = 1 - \sum_{k \geq n} \pi[\mathbb{N}+j, a](k).$$

Suppose for a contradiction that $x_0 > b$. It is then impossible to find $x_1, x_2, \ldots$ in $[0, 1]$ such that
$$\sum_{k \geq 0} x_k = 1, \qquad \sum_{k \geq 0} k x_k = R \quad \text{and} \quad (x_0, x_1, x_2, \ldots) \leq_{\text{st}} \pi[\mathbb{N}+j, a].$$

Therefore $x_0 \leq b$. We can take $n = \mathrm{O}(R \log R)$, so the result follows. □

**9. Convergence of $D_{\mathrm{KL}} = D_{\mathrm{KL}}(\mathbf{x} \| G^R)$.** Now that we have bounds on the fluid limit, we can show that $D_{\mathrm{KL}}$ decreases rapidly. As soon as we have the lower bound $x_0 \geq 1/(5R)$, we can use Lemma 7.1.

LEMMA 9.1. *As $t \to \infty$,*
$$D_{\mathrm{KL}}(t) = D_{\mathrm{KL}}(0) \exp(-\Omega(R^{-1}\sqrt{t})), \qquad D_{\mathrm{KL}}(0) = \mathrm{O}(\log R).$$

This eliminates the possibility that $D_{\mathrm{KL}}$ gets "stuck" some distance away from zero. Once $D_{\mathrm{KL}}$ is sufficiently small, we can use Lemma 7.2 to prove Theorem 1.2.

PROOF OF LEMMA 9.1. As discussed in Section 7, the log Sobolev constant of BRW$[n, a]$, $\alpha = \Theta(a/n)$. If $(1 - x_0)^n < 1/2$, by the first line of (5.3), (5.4) and Lemma 7.1 with $\mu = \mathbf{x}$ and $a = x_0$,
$$\frac{d}{dt}S \geq \Theta(x_0/n)\left[D_{\mathrm{KL}} - \sum_{k \geq n} x_k \log \frac{x_k}{x_0(1-x_0)^k} - 2(1-x_0)^n\right].$$

If we can find $n = n(D_{\mathrm{KL}})$ such that, say,

(9.1) $$2(1-x_0)^n + \sum_{k \geq n} x_k \log \frac{x_k}{x_0(1-x_0)^k} \leq D_{\mathrm{KL}}/2,$$

then

(9.2) $$\frac{d}{dt}S = \frac{d}{dt}(-D_{\mathrm{KL}}) \geq \Theta(x_0/n)D_{\mathrm{KL}}.$$

We must let $n \to \infty$ as $D_{\mathrm{KL}} \to 0$ to maintain inequality (9.1). However, the quicker $n$ grows, the weaker the bound (9.2) becomes.

To apply Lemma 8.3, assume $t \geq s_0$; let $a = 1/(5R)$, $b = 1 - (R+1)^{-O(1)}$ and $j = O(R \log R)$. We then have $x_0 \in [a, b]$ and $\mathbf{x} \leq_{\mathrm{st}} \pi[\mathbb{N} + j, a]$;

$$2(1-x_0)^n + \sum_{k \geq n} x_k \log \frac{x_k}{x_0(1-x_0)^k}$$

$$\leq 2(1-a)^n + \sum_{k \geq n} \pi[\mathbb{N}+j, a](k)[-\log a - k \log(1-b)]$$

$$\leq 2(1-a)^n + (1-a)^{n-j}[-\log a - (n+a^{-1})\log(1-b)].$$

Of course, $1 - a \leq \exp(-a)$. Taking $n = O(-R \log D_{\mathrm{KL}})$ above, we satisfy inequality (9.1). Putting $n = O(-R \log D_{\mathrm{KL}})$ and $x_0 = \Omega(R^{-1})$ into (9.2),

$$\frac{d}{dt}(-\log D_{\mathrm{KL}}) \geq \Theta(R^{-2}/(-\log D_{\mathrm{KL}})). \qquad \square$$

PROOF OF THEOREM 1.2.　We start off by preparing a "daisy-chain" of inequalities, using (5.3), (5.4), Lemma 7.2 and the function $\phi$:

(i) $\frac{d}{dt}S = \frac{d}{dt}(-D_{\mathrm{KL}}) \geq \sum_{k=0}^{n-2} \frac{(x_k(1-x_0)-x_{k+1})^2}{\max\{x_k(1-x_0), x_{k+1}\}}$,

(ii) $\sum_{k=0}^{n-2} \frac{(x_k(1-x_0)-x_{k+1})^2}{\max\{x_k(1-x_0), x_{k+1}\}} \geq \sum_{k=0}^{n-2} \frac{(x_k(1-x_0)-x_{k+1})^2}{x_0(1-x_0)^{k+1}} \inf_{k<n} \frac{x_0(1-x_0)^k}{x_k}$,

(iii) $\sum_{k=0}^{n-2} \frac{(x_k(1-x_0)-x_{k+1})^2}{x_0(1-x_0)^{k+1}} \geq (x_0^2/4) \sum_{k=1}^{n-1} \frac{(x_k - x_0(1-x_0)^k)^2}{x_0(1-x_0)^k}$,

(iv) $\sum_{k=1}^{n-1} \frac{(x_k - x_0(1-x_0)^k)^2}{x_0(1-x_0)^k} \geq \sum_{k=1}^{n-1} \phi(x_0(1-x_0)^k, x_k)$,

(v) $\sum_{k=1}^{\infty} \phi(x_0(1-x_0)^k, x_k) \geq D_{\mathrm{KL}}$.

Joining these together,

$$\frac{d}{dt}(-D_{\mathrm{KL}}) \geq \frac{x_0^2}{4}\left(D_{\mathrm{KL}} - \sum_{k \geq n} \phi(x_0(1-x_0)^k, x_k)\right) \inf_{k<n} \frac{x_0(1-x_0)^k}{x_k}.$$

The result now follows if we can choose $n = n(D_{\mathrm{KL}})$ such that as $D_{\mathrm{KL}} \to 0$,

(I) $\displaystyle \frac{1}{D_{\mathrm{KL}}} \sum_{k \geq n} \phi(x_0(1-x_0)^k, x_k) \to 0$, 　　(II) $\displaystyle \inf_{k \leq n} \frac{x_0(1-x_0)^k}{x_k} \geq 1/2$.

We need $n$ large for (I), but not too large or (II) might fail.

For $\ell = 1, 2, \ldots$, let $c_\ell = \exp(-\exp(\ell))$. We will use the sequence $(c_\ell)$ to measure the rate of decrease of $D_{\mathrm{KL}}$. As $D_{\mathrm{KL}}$ is monotonically decreasing, we can let $s_\ell$ be the unique time at which $D_{\mathrm{KL}} = c_\ell$. We have chosen $(c_\ell)$ to be decreasing doubly exponentially, so we must show that $(s_\ell)$ grows only exponentially.

Consider the time period $s_\ell$ to $s_{\ell+1}$. Using inequality (8.3) with $D_{\text{KL}} \leq c_\ell$, we can find $a$ and $b$ such that $x_0 \in [a, b]$ from $s_\ell$ onward:

$$(9.3) \qquad a = \frac{1}{R+1} - \sqrt{2c_\ell}, \qquad b = \frac{1}{R+1} + \sqrt{2c_\ell}.$$

Let $\mu_a$ and $\mu_b$ be defined as in Lemma 8.2 with $s = s_\ell$. To take advantage of the bounds $\mu_b \leq_{\text{st}} \mathbf{x} \leq_{\text{st}} \mu_a$, we must wait until $\mu_a$ and $\mu_b$ are close to their respective equilibrium distributions. We will split the time period $[s_\ell, s_{\ell+1}]$ into two parts: $[s_\ell, s'_\ell]$ where we wait for $\mu_a$ and $\mu_b$ to mix and $[s'_\ell, s_{\ell+1}]$, where (I) and (II) hold with $n = \lceil \log(c_{\ell+2})/\log(1-a) \rceil$.

By inequality (8.1), from time $s'_\ell := s_\ell + \mathrm{O}(R^2 e^\ell)$ onward, the bound $\mu_a$ satisfies (and similarly for $\mu_b$)

$$(9.4) \qquad \|\mu_a - \pi[\mathbb{N}, a]\|_{\text{TV}} \leq \|\mu_a - \pi[\mathbb{N}, a]\|_{(1)} \leq c_\ell^2.$$

Combining (9.4) with $\mathbf{x} \leq_{\text{st}} \mu_a$ allows us to bound the tail of $\mathbf{x}$,

$$\sum_{k \geq n} k x_k \leq \sum_{k \geq n} k \mu(k)$$
$$\leq \sum_{k \geq n} k \pi[\mathbb{N}, a](k) + \|\mu_a - \pi[\mathbb{N}, a]\|_{(1)}$$
$$\leq (n + a^{-1})(1-a)^n + c_\ell^2.$$

Let $t \in [s'_\ell, s_{\ell+1}]$. By the definition of $\phi$,

$$\sum_{k \geq n} \phi(x_0(1-x_0)^k, x_k) = \sum_{k \geq n} x_k \log \frac{x_k}{x_0(1-x_0)^k} - (x_k - x_0(1-x_0)^k)$$
$$\leq (1-a)^n + \sum_{k \geq n} x_k[-\log a - k\log(1-b)].$$

As $(1-a)^n \leq c_{\ell+2} \ll c_{\ell+1} \leq D_{\text{KL}}$, this implies (I). By inequality (8.2),

$$\frac{x_0(1-x_0)^k}{x_k} \geq \frac{a(1-b)^k}{(1-a)^k - (1-b)^{k+1} + 2c_\ell^2}.$$

Implicitly, $n$ is a function of $\ell$. By (9.3), in limit as $\ell \to \infty$,

$$\sup_{k < n}(1-a)^k/(1-b)^k \to 1$$

so (II) holds. For $\ell$ sufficiently large

$$\frac{d}{dt}(-\log D_{\text{KL}}) = \Omega(R^{-2}).$$

This gives $s_{\ell+1} - s'_\ell = \mathrm{O}(R^2 e^\ell)$ and therefore $s_{\ell+1} - s_\ell = \mathrm{O}(R^2 e^\ell)$, too. Hence $s_\ell = \mathrm{O}(R^2 e^\ell)$. $\square$

**10. Convergence to the fluid limit.** We now prove that the empirical distribution converges to the fluid limit.

PROOF OF THEOREM 1.1. We will use Theorem 2.2 from [7]. For $N \in \mathbb{N}$, the empirical distribution is a pure jump process in $I_N = (\mathbb{N}/N)^{\mathbb{N}}$. Let vectors $\mathbf{e}_0 = (1, 0, 0, \ldots)$, $\mathbf{e}_1 = (0, 1, 0, \ldots), \ldots$ be the canonical basis for $\mathbb{R}^{\mathbb{N}}$. The Lévy kernel of $\mathbf{X}$ is defined for distributions $\mathbf{x} \in I_N$ with mass $\sum x_k = 1$,

$$K^N(\mathbf{x}, d\mathbf{y}) = N \sum_{i>0} \sum_{j \geq 0} \left( x_i x_j - \frac{x_i 1_{\{i=j\}}}{N} \right) \delta_{[\mathbf{e}_{i-1}/N - \mathbf{e}_i/N - \mathbf{e}_j/N + \mathbf{e}_{j+1}/N]}.$$

The formal "limit kernel" is

$$K(\mathbf{x}, d\mathbf{y}) = \lim_{N \to \infty} \frac{K^N(\mathbf{x}, d\mathbf{y}/N)}{N} = \sum_{i>0} \sum_{j \geq 0} x_i x_j \delta_{[\mathbf{e}_{i-1} - \mathbf{e}_i - \mathbf{e}_j + \mathbf{e}_{j+1}]}.$$

The limit kernel encodes the fluid limit differential equation (5.2),

$$\frac{d}{dt}\mathbf{x}(t) = b(\mathbf{x}(t)), \qquad b(\mathbf{x}) := \int \mathbf{y} K(\mathbf{x}, d\mathbf{y}).$$

We must first show that the fluid limit differential equation has a unique solution. Consider the set $S = \{\mathbf{x} \in [0,1]^{\mathbb{N}} : \sum_k x_k \leq 1\}$ of "subprobability" distributions. Equipped with the sup-norm, $S$ is a complete metric space and $b$ is a Lipschitz function on $S$.

Now fix $T > 0$. To apply the fluid limit theorem, we must bound the tails of the fluid limit and the empirical distribution for $t \in [0, T]$. Let $\mu_t$ be the distribution of a random walk on $\mathbb{N}$ that starts at $R$, then steps forward at rate 1. This is a Poisson arrival process; the number of steps taken forward is Poisson($t$). By a slight modification to Lemma 8.2,

$$\mathbf{x}(t) \leq_{\text{st}} \mu_t \leq_{\text{st}} \mu_T.$$

The empirical distribution has an analogous bound. Recall that $B_1, \ldots, B_N$ are the numbers of balls in each box. Let $C_1, \ldots, C_N$ be $N$ independent copies of the Poisson arrival process, each starting $C_i(0) = R$. There is a coupling $(\mathbf{B}, \mathbf{C})$ such that $B_i \leq C_i$ for all $i$: every time box $i$ is picked as the sink box in the MFZRP, increment $C_i$ by one.

Let $d \in \mathbb{N}$, $d > R$. Let $p(d)$ be the probability that $C_i(T) \geq d$. The number of boxes with occupancy exceeding $d$ at any point up to time $T$ is stochastically smaller than binomial $\text{Bin}(N, p(d))$. By Markov's inequality, $p(d) \to 0$ as $d \to \infty$. By concentration, Theorem 8.6, the probability that $\sum_{k \geq d} X_k$ exceeds $2p_d$ at any point up to time $T$ decays exponentially in $N$.

The proof of Theorem 2.2 [7] considers truncated, finite-dimensional versions of the jump process and fluid limit. Truncating the distribution (setting $x_k = 0$ for $k \geq d$) changes the mass to $m = \sum x_k \in [0, 1]$. We can use the formula for $K^N$ to extend the MFZRP to all points $\mathbf{x} \in I_N \cap S$. The corresponding balls and boxes Markov chain has $Nm$ boxes, with balls distributed according to $\mathbf{x}$. Each ordered pair of boxes interacts at rate $N^{-1}$. If $m = 1$, the process is unchanged. □

**11. The initial configuration and Markov chain mixing.** We have studied the relaxation of the MFZRP by showing that it tends to a fluid limit as the number of boxes tends to infinity, and that the fluid limit converges exponentially to $G^R$. Combining Theorems 1.1 and 1.2, we can conclude that with high probability as $N \to \infty$, the empirical distribution rapidly approaches the geometric distribution $G^R$.

An alternative, and in a sense stronger description would be provided by finding the total variation mixing time for the finite MFZRP. The process can only be close to the equilibrium distribution on $\mathcal{B}_N$ in total variation if the empirical distribution is close to $G^R$ on $\mathbb{N}$ in Kullback–Leibler divergence. This follows from (3.1)–(3.3): if with probability close to 1, $D_{\mathrm{KL}}(\mathbf{X}(t) \| G^R)$ is significantly larger than expected equilibrium value, the process cannot be close to equilibrium at time $t$.

When studying the fluid limit we restricted our attention to the initial distribution that arises when every box starts with $R$ balls in: $\mathbf{x}(0) = \delta_R$. This suggests we ask two different questions about the finite MFZRP. First, what is the (uniform) total variation mixing time $\tau_1(1/4)$, and second, what is the total variation mixing time when started with $\mathbf{X}(0) = \delta_R$?

The answer to the first question is order between $NR$ and $NR^2 \log R$. In [16], the spectral gap for the MFZRP is shown to have order $\Omega(R^{-2})$ uniformly in $N$. Combined with (6.3) and $-\log \pi_{\min} = \mathrm{O}(N \log R)$, this provides the upper bound.

Under the MFZRP equilibrium measure, the number of balls in each box is $\mathrm{O}(R \log N)$ with probability tending to 1 as $N \to \infty$. This provides a simple lower bound on the convergence time. Start all $NR$ balls in box 1, with the other $N-1$ boxes empty. While with high probability box 1 contains more than $\mathrm{O}(R \log N)$ balls, the process cannot be close to equilibrium. Box 1 is losing balls at rate 1, so time $\Omega(NR)$ is needed to get rid of the excess $NR - \mathrm{O}(R \log N)$ balls.

We do not have a satisfactory answer to the second question. In [8], rather than taking the fluid limit, we worked with the entropy of the empirical distribution directly. With probability tending to 1 as $N \to \infty$, $D_{\mathrm{KL}}(\mathbf{X}(t) \| G^R)$ decreases to $\mathrm{O}(R^4 N^{-1} \log N)$ in time $\mathrm{O}(R^3 \log N)$. With a little extra care, both these bounds can be improved by a factor of $R$.

**Acknowledgments.** The bulk of this work was done as part of my Ph.D. [8], Chapter 4. I would like to thank my Ph.D. supervisor Geoffrey Grimmett and Archie Howie for suggesting the problem.

I would also like to thank my Ph.D. examiners, James Martin and Yuri Suhov, for their helpful comments on my thesis, and the referee for suggesting a number of ways of simplifying the proofs.

DEPARTMENT OF MATHEMATICS
UNIVERSITY OF BRITISH COLUMBIA
1984 MATHEMATICS ROAD
VANCOUVER, BRITISH COLUMBIA
CANADA V6T 1Z2
E-MAIL: ben@math.ubc.ca